\documentclass[12pt]{amsproc}
\usepackage[foot]{amsaddr}
\usepackage{amsmath,amssymb,amsthm}
\usepackage{color,graphics,srcltx,float}
\usepackage{placeins}
\usepackage[left=1in,right=1in,top=1in,bottom=1in]{geometry}
\usepackage{stmaryrd}
\usepackage{enumerate}
\usepackage{booktabs}
\usepackage[table,x11names]{xcolor}
\usepackage{longtable}
\usepackage{makecell}
\usepackage{hyperref}
\usepackage{soul}
\hypersetup{
    colorlinks,
    linkcolor={red!50!black},
    citecolor={blue!50!black},
    urlcolor={blue!80!black}
}
\usepackage{algpseudocode}

\definecolor{darkblue}{rgb}{0,0,0.8}

\newtheorem{theorem}{Theorem}[section]

\newtheorem{corollary}[theorem]{Corollary}

\theoremstyle{definition}

\renewcommand{\leq}{\leqslant}

\usepackage{color}

\allowdisplaybreaks

\title{A census of small Schurian association schemes}
\author[Lansdown]{Jesse Lansdown}

\address{School of Mathematics and Statistics, University of Canterbury, Christchurch, New Zealand}
\address{Centre for the Mathematics of Symmetry and Computation, Department of Mathematics and Statistics, The University of Western Australia, Crawley, WA 6009, Australia.}
\email{jesse.lansdown@canterbury.ac.nz}

\thanks{The author was supported by the Australian Research Council Discovery Grant DP200101951. This work was supported by resources provided by the Pawsey Supercomputing Centre with funding from the Australian Government and the Government of Western Australia.}

\begin{document}

\maketitle

\begin{abstract}
Using the classification of transitive groups of degree $n$, for $2 \leqslant n \leqslant 48$,
we classify the Schurian association schemes of order $n$, and as a consequence, the transitive groups of degree $n$ that are $2$-closed. In addition, we compute the character table of each association scheme and provide a census of important properties. Finally, we compute the $2$-closure of each transitive group of degree $n$, for $2 \leqslant n \leqslant 48$. The results of this classification are made available as a supplementary database.
\end{abstract}

\section{Introduction}

Association schemes are some of the most important objects in algebraic combinatorics, with applications to coding theory, finite geometry, group theory, and even statistics.
They generalise the concepts of strongly-regular and distance-regular graphs to describe structures with high degrees of regularity. As a result they arise naturally in many settings.
Moreover, Delsarte showed that many combinatorial objects can be described by the eigenspaces of an association scheme \cite{delsarte_algebraic_1973}, providing a powerful tool for studying geometric objects, cliques of graphs, designs, codes, and more.

The (not necessarily commutative) association schemes of order $n$ have been classified for $1 \leqslant n \leqslant 34$ and also for $n = 38$ by Nomiyama \cite{Assoc10}, Hirasaka \cite{Assoc11and12}, Hirasaka and Suga \cite{Assoc13and15}, 
Hanaki and Miyamoto \cite{Assoc16and17,Assoc18and19,AssocSmall}, and Hanaki, Kharaghani, Mohammadian, and Tayfeh-Rezaie \cite{Assoc31}. With the exception of very small values of $n$, this classification necessitated the use of a computer and further classification is extremely difficult due to combinatorial explosion.
Hanaki maintains a website with these association schemes at \cite{hanaki_miyamoto}. This database has been used for forming and testing conjectures on association schemes.

The regularity properties of association schemes may be thought of as capturing combinatorial symmetry.
Indeed, their connection to group actions was a crucial motivator in their development \cite{Higman1987}.
Association schemes with the strongest connection to groups are called \emph{Schurian}; their relations correspond to the orbitals of a transitive group. 
Schurian association schemes remain one of the strongest tools for studying group actions and combinatorial objects with high degrees of symmetry. For example, they have found recent use in exploring the synchronisation hierarchy of permutation groups (e.g. \cite{Diagonal}).

The transitive groups of degree $n$ have been classified for $1 \leq n \leq 48$ by Miller \cite{MillerCollected,Miller1896}, Royle \cite{Royle1987}, Hulpke \cite{Hulpke2005}, Cannon and Holt \cite{CannonHolt2008}, Holt and Royle \cite{TrnsGrp1to47}, and Holt, Royle, and Tracey \cite{TrnsGrp48}, mostly by computer, with the exception of some small values of $n$.
For $n=48$ alone, there are $195,826,352$ conjugacy classes of transitive groups, making their classification a particular milestone in the classification of transitive permutation groups of low degree.
This recent work on classifying transitive groups makes the main result of this paper possible.

\begin{theorem}\label{thm:main}
There are are $24678$ Schurian association schemes of order $2$ to $48$ up to isomorphism.
\end{theorem}

There is a bijection between the Schurian association schemes and the $2$-closed transitive groups, which yields the following:

\begin{corollary}
There are $24678$  $2$-closed groups of degree $2$ to $48$ up to isomorphism.
\end{corollary}

We provide the association schemes corresponding to Theorem \ref{thm:main} as a supplement \cite{SmallSchurianSchemesDatabase} to this paper.
It is hoped that this database will prove as useful for studying Schurian association schemes and $2$-closed groups as \cite{hanaki_miyamoto} has for general association schemes.
We also compute the character table of each association scheme since this is a computationally intensive, and potentially prohibitive, task.
A census of some of the important properties is provided in Table \ref{table:AllValues}, where the total number of association schemes of each order is given, along with the number which are stratifiable, commutative, symmetric, primitive, metric, cometric, and thin, respectively.
The rows that are highlighted indicate previously unclassified Schurian association schemes, which are not available in \cite{hanaki_miyamoto} or elsewhere.

We also identify and provide as an additional supplement \cite{TwoClosuresDatabase} containing the $2$-closures for all transitive groups of degree $2$ to $48$. This supplementary data is likely to be of special interest in the case $n=48$ since transitive identification is not computable in \emph{Magma} \cite{Magma}  or \emph{GAP} \cite{GAP4} for this degree.

The computations required by the classification are described in Section \ref{Computation}. In particular, the fast practical computation of automorphisms and isomorphisms is made possible by representing association schemes as suitable digraphs. Hence the $2$-closures of the corresponding groups may also be computed quickly, which is otherwise difficult and slow.

\section{Association schemes} \label{Background}

The term \emph{association scheme} is used in different ways throughout the literature\footnote{Peter Cameron discusses the differing ``association scheme'' terminology on his blog:\\ \href{https://cameroncounts.wordpress.com/2014/06/08/terminology-association-scheme-or-coherent-configuration/}{https://cameroncounts.wordpress.com/2014/06/08/terminology-association-scheme-or-coherent-configuration/}}. The definition used in this paper is also known as a \emph{homogeneous coherent configuration}. We use this definition because it is the most general and is used for the corresponding objects in the catalogue of association schemes at \cite{hanaki_miyamoto}. We refer to \cite{BannaiIto} for greater detail on association schemes, but provide some of the relevant definitions and results (without proof) in this section. The connection between association schemes and permutation groups is explored in greater detail in \cite{CAP}.

Let $\Omega$ be a finite set of cardinality $n$ and let $\mathcal{R} = \{ R_0, R_1, \ldots, R_d\}$ be subsets of $\Omega \times \Omega$. We shall call $R_i \in \mathcal{R}$ a \emph{relation} and refer to $R_{i}^\top = \{ (y, x) : (x, y) \in R_i \}$ as its \emph{converse relation}.
Then $(\Omega, \mathcal{R})$ is an \emph{association scheme} with $d$ \emph{classes} and \emph{order} $n$ if the following hold:
\begin{enumerate}
	\item $\mathcal{R}$ is a partition of $\Omega \times \Omega$,
	\item $R_0 = \{ (x, x) : x \in \Omega \}$,
	\item $R_i^\top \in \mathcal{R}$ for all $R_i \in \mathcal{R}$,
	\item there exist constants $p_{ij}^k$ (called \emph{intersection numbers}) for $0 \leq i, j, k \leq d$ such that for any $(x, y) \in R_k$, 
	\[
		p_{ij}^k = | \{ z \in \Omega : (x, z) \in R_i, (z, y) \in R_j \} |.
	\]
\end{enumerate}

Moreover an association scheme is \emph{symmetric} if $R_i^\top = R_i$ for all $0 \leq i \leq d$, \emph{commutative} if $p_{ij}^k = p_{ji}^k$ for all $0 \leq i, j, k \leq d$, and \emph{stratifiable} 
if an association scheme is formed by replacing each relation and its converse by their union.
A scheme is called \emph{thin} if $d = n-1$, in which case it is equivalent to a group.
Two association schemes $(\Omega, \mathcal{R})$ and $(\Delta, \mathcal{S})$ are \emph{isomorphic} if they have the same number of classes, $d$, and there exist bijections $f:\Omega \to \Delta$ and $g: \{0, \ldots d\} \to \{0, \ldots d\}$ such that $(x, y) \in R_i$ if and only if $(f(x), f(y)) \in S_{g(i)}$.

Let $G$ be a transitive permutation group acting on $\Omega$. Then the orbits of $G$ on $\Omega \times \Omega$ form the relations of an association scheme, $\mathcal{K}(G)$. An association scheme $(\Omega, \mathcal{R})$ is called \emph{Schurian} if $(\Omega, \mathcal{R}) = \mathcal{K}(G)$ for some transitive permutation group $G$.
An \emph{automorphism} of an association scheme $(\Omega, \mathcal{R})$ is a permutation of $\Omega$ which fixes every $R_i \in \mathcal{R}$. The set of all automorphisms forms the \emph{automorphism group}, ${\rm Aut}((\Omega, \mathcal{R}))$. 
Note that $G \leqslant {\rm Aut}(\mathcal{K}(G))$. Moreover, $G$ is $2$-closed in the case of equality\footnote{The $2$-closure of a group $G$ is the largest subgroup of ${\rm Sym}(\Omega)$ preserving the orbits of $G$ on $\Omega \times \Omega$ and $G$ is called $2$-closed if it equals its $2$-closure.}.

The \emph{adjacency matrix} with respect to $R_i$ is the $n \times n$ matrix $A_i$ indexed by the elements of $\Omega$, where
\[ (A_i)_{xy} = \begin{cases} 
      1, &  \text{if } (x, y) \in R_i; \\
      0, &  \text{otherwise}.
   \end{cases}
\]
If the digraph defined by the adjacency matrix $A_i$ is connected for all $1 \leq i \leq d$, then we call $(\Omega, \mathcal{R})$ \emph{primitive}.
To store an association scheme compactly, we define the \emph{relation matrix},
\[
M := \sum_{i=0}^d  i A_i.
\]
Note that since the relations partition $\Omega \times \Omega$, we can recover each $A_i$ and hence the relations $\mathcal{R}$ for $(\Omega,\mathcal{R})$ from the relation matrix $M$.

For an association scheme $(\Omega, \mathcal{R})$, the adjacency matrices span a semisimple $\mathbb{C}$-algebra called the \emph{adjacency algebra} or the \emph{Bose-Mesner} algebra, denoted $\mathbb{C}\mathcal{R}$. When $(\Omega, \mathcal{R})$ is commutative, there exists a second basis $\{E_0, \ldots, E_d\}$ for $\mathbb{C}\mathcal{R}$ consisting of minimal idempotents (ie. $E_i E_j = 0$ for $i \neq j$, and $E_i E_i = E_i$).
A \emph{representation} $\varphi$ of an association scheme $(\Omega, \mathcal{R})$ is an algebra homomorphism from the adjacency algebra $\mathbb{C}\mathcal{R}$ to the full matrix algebra over $\mathbb{C}$ and we define the \emph{character} $\chi$ afforded by $\varphi$ by $\chi(A_i) = {\rm trace}(\varphi(A_i))$. 
The set of irreducible characters of $\mathbb{C}\mathcal{R}$ is denoted by ${\rm Irr}(\mathcal{R})$. The \emph{standard representation $\Gamma_{\mathcal{R}}$} is the representation which sends each adjacency matrix to itself, that is, $\Gamma_\mathcal{R}(A_i) = A_i$ for all $0 \leqslant i \leqslant d$.
The \emph{standard character} $\gamma_\mathcal{R}$ is the character afforded by the standard representation $\Gamma_\mathcal{R}$ and satisfies
\[
\gamma_\mathcal{R}(A_i) = 
\begin{cases} 
      |\Omega| 1_\mathbb{C}, &   i = 0; \\
      0, &  \text{otherwise},
   \end{cases}
\]
and has the irreducible decomposition
\[
\gamma_\mathcal{R} = \sum_{\chi \in {\rm Irr}(\mathcal{R})} m_\chi . \chi,
\]
where $m_\chi$ is called the multiplicity of $\chi$. Let $T$ be the $|{\rm Irr}(\mathcal{R})| \times (d+1)$ matrix with rows indexed by ${\rm Irr}(\mathcal{R})$ and columns by $\mathcal{R}$ such that the $T_{\chi, R_i} = \chi (A_i)$. Then $T$ is called the \emph{character table} of $(\Omega, \mathcal{R})$.

Finally, let $(\Omega,\mathcal{R})$ be symmetric.
If there is an ordering $\{A_i\}_{i=0}^d$ such that there exist polynomials $v_i$ of degree $i$ for $i \in \{0, \ldots, d\}$ with the property $v_i(A_1)=A_i$, then we call $(\Omega, \mathcal{R})$ \emph{P-polynomial} or \emph{metric}. Similarly, if there is an ordering $\{E_i\}_{i=0}^d$ such that there exist polynomials $q_i$ of degree $i$ for $i \in \{0, \ldots, d\}$ with the property $q_i(E_1)=E_i$, then we call $(\Omega, \mathcal{R})$ \emph{Q-polynomial} or \emph{cometric}. Metric and cometric association schemes have interesting algebraic and combinatorial properties and are often studied in their own right.

\section{Computation} \label{Computation}

By definition, an association scheme $(\Omega, \mathcal{R})$ is Schurian precisely when there exists a group $G$ acting transitively on $\Omega$ such that $(\Omega, \mathcal{R}) = \mathcal{K}(G)$.
Let $\mathcal{G}$ be a set of transitive permutation groups of degree $n$ such that every conjugacy class of transitive subgroups of $S_n$ has precisely one representative in $\mathcal{G}$.
Then $G$ is conjugate in $S_n$ to some group in $\mathcal{G}$, and without loss of generality we may assume that $G \in \mathcal{G}$.

Note that for two non-isomorphic groups $G_1$ and $G_2$ in $\mathcal{G}$ it is possible that $\mathcal{K}(G_1) \cong \mathcal{K}(G_2)$. However, we may use the fact that $G \leq {\rm Aut}(\mathcal{K}(G))$ to determine non-isomorphic representatives of all Schurian association schemes of order $n$ by only keeping $\mathcal{K}(G)$ when $G$ is the full automorphism group. Hence we derive the classification of Schurian association schemes of order $n$ from the classification of transitive permutation groups of degree $n$ by the following procedure:
\vspace{0.5cm}

\begin{algorithmic}
\State $L \gets \{\}$
\For{$G \in \mathcal{G}$}
	\If{$G = {\rm Aut}(\mathcal{K}(G))$} 
		\State $L \gets L \cup \{\mathcal{K}(G)\}$
	\EndIf
\EndFor
\end{algorithmic}
\vspace{0.5cm}

The transitive permutation groups of degree at most $48$ are available in both \emph{MAGMA} \cite{Magma} and \emph{GAP} \cite{GAP4} (via the \emph{TransGrp} package \cite{TransGrp3.6}). In both cases, the groups of degree $32$ and $48$ must be downloaded\footnote{Transitive groups of degree $32$ and $48$ in \emph{Magma}: \url{http://magma.maths.usyd.edu.au/magma/download/db/}, Transitive groups of degree $32$ in \emph{GAP}: \url{https://www.math.colostate.edu/~hulpke/transgrp/trans32.tgz}, Transitive groups of degree $48$ in \emph{GAP}: \url{https://zenodo.org/record/5935751}} separately due to their size. The computations for this paper were done in \emph{GAP}, making use of the author's \emph{GAP} package, \emph{AssociationSchemes} \cite{AssociationSchemes}, which among other things can construct $\mathcal{K}(G)$, find the automorphism group of an association scheme, and compute isomorphisms between association schemes.

The $2$-closure of a group $G$ is precisely ${\rm Aut}(\mathcal{K}(G))$, and so classifying the Schurian association schemes is equivalent to classifying the $2$-closed transitive permutation groups. In \emph{GAP} the $2$-closure of a group can be computed  using the \emph{GAP} package \emph{GRAPE} \cite{GRAPE}, however this is significantly slower than computing the automorphism group of $\mathcal{K}(G)$ using \emph{AssociationSchemes} directly, although it does provide a means of verification. In practice, it may be faster to check that ${\rm Aut}(\mathcal{K}(G)) \leq G$ by testing containment of generators of ${\rm Aut}(\mathcal{K}(G))$ in $G$. To accelerate the computation, the groups in $\mathcal{G}$ were divided and multiple jobs were run in parallel. %\textcolor{blue}{mention Nimbus?}
{\rm 
The character table and other properties of the Schurian association schemes were also computed using the \emph{AssociationSchemes} package. In fact, some improvements have been made to the package in the process. However, manual intervention was required in many difficult cases.

The number of association schemes of each type with a given order are provided in Table \ref{table:AllValues}. The number of Schurian association schemes which are \emph{stratifiable}, \emph{commutative}, \emph{symmetric}, \emph{primitive}, \emph{metric}, \emph{cometric}, and \emph{non-Schurian} are abbreviated as \emph{strat.}, \emph{com.}, \emph{sym.}, \emph{prim.}, \emph{met.},  \emph{comet.}, and \emph{NS} respectively. 
{\small
\begin{table}[htbp]
\centering
%Order & Total & Strat. & Comm. & Sym. & Primitive & Metric & Cometric & Thin \\
\begin{tabular}{ r | r r r r r r r r | r }
\toprule
 Order & Total & Strat. & Com. & Sym. & Prim. & Met. & Comet. & Thin & NS\\
\midrule
2 &1 & 1 & 1 & 1 & 1 & 1 & 1 & 1 & 0 \\
3 & 2 & 2 & 2 & 1 & 2 & 1 & 1 & 1 & 0 \\
4 & 4 & 4 & 4 & 3 & 1 & 2 & 2 & 2 & 0 \\
5 & 3 & 3 & 3 & 2 & 3 & 2 & 2 & 1 & 0 \\
6 & 8 & 7 & 7 & 4 & 1 & 4 & 4 & 2 & 0 \\
7 & 4 & 4 & 4 & 2 & 4 & 2 & 2 & 1 & 0 \\
8 & 21 & 20 & 19 & 10 & 1 & 5 & 5 & 5 & 0 \\
9 & 12 & 12 & 12 & 6 & 2 & 4 & 4 & 2 & 0 \\
10 & 13 & 11 & 11 & 8 & 2 & 6 & 6 & 2 & 0 \\
11 & 4 & 4 & 4 & 2 & 4 & 2 & 2 & 1 & 0 \\
12 & 59 & 47 & 47 & 21 & 1 & 8 & 8 & 5 & 0 \\
13 & 6 & 6 & 6 & 4 & 6 & 3 & 3 & 1 & 0 \\
14 & 16 & 14 & 14 & 8 & 1 & 6 & 6 & 2 & 0 \\
15 & 24 & 23 & 23 & 10 & 2 & 6 & 5 & 1 & 1 \\
16 & 206 & 171 & 158 & 56 & 4 & 9 & 9 & 14 & 16 \\
17 & 5 & 5 & 5 & 4 & 5 & 3 & 3 & 1 & 0 \\
18 & 93 & 71 & 71 & 32 & 1 & 8 & 7 & 5 & 2 \\
19 & 6 & 6 & 6 & 3 & 6 & 2 & 2 & 1 & 1 \\
20 & 95 & 73 & 73 & 41 & 1 & 10 & 8 & 5 & 0 \\
21 & 32 & 29 & 29 & 11 & 3 & 6 & 5 & 2 & 0 \\
22 & 16 & 14 & 14 & 8 & 1 & 6 & 6 & 2 & 0 \\
23 & 4 & 4 & 4 & 2 & 4 & 2 & 2 & 1 & 18 \\
24 & 669 & 454 & 438 & 136 & 1 & 9 & 10 & 15 & 81 \\
25 & 32 & 32 & 32 & 20 & 9 & 5 & 5 & 2 & 13 \\
26 & 24 & 20 & 20 & 14 & 1 & 6 & 6 & 2 & 10 \\
27 & 122 & 112 & 112 & 38 & 5 & 7 & 6 & 5 & 380 \\
28 & 124 & 103 & 103 & 47 & 4 & 10 & 9 & 4 & 61 \\
29 & 6 & 6 & 6 & 4 & 6 & 3 & 3 & 1 & 20 \\
30 & 228 & 166 & 166 & 73 & 1 & 11 & 10 & 4 & 15 \\
31 & 8 & 8 & 8 & 4 & 8 & 2 & 2 & 1 & 98299 \\
32 & 4261 & 2579 & 2264 & 413 & 1 & 13 & 11 & 51 & 13949 \\
33 & 27 & 27 & 27 & 9 & 1 & 4 & 4 & 1 & 0 \\
34 & 20 & 16 & 16 & 13 & 1 & 5 & 5 & 2 & 0 \\
\rowcolor{black!80!brown!15} 35 & 43 & 43 & 43 & 17 & 3 & 6 & 6 & 1 &  \\
\rowcolor{black!80!brown!15} 36 & 1274 & 806 & 804 & 276 & 9 & 17 & 15 & 14  & \\
\rowcolor{black!80!brown!15} 37 & 9 & 9 & 9 & 6 & 9 & 3 & 3 & 1  & \\
38 & 22 & 19 & 19 & 10 & 1 & 5 & 5 & 2  & 11 \\
\rowcolor{black!80!brown!15} 39 & 44 & 41 & 41 & 15 & 1 & 4 & 4 & 2  & \\
\rowcolor{black!80!brown!15} 40 & 1095 & 712 & 687 & 262 & 3 & 11 & 11 & 14  & \\
\rowcolor{black!80!brown!15} 41 & 8 & 8 & 8 & 6 & 8 & 3 & 3 & 1  & \\
\rowcolor{black!80!brown!15} 42 & 298 & 210 & 210 & 81 & 1 & 11 & 10 & 6  & \\
\rowcolor{black!80!brown!15} 43 & 8 & 8 & 8 & 4 & 8 & 2 & 2 & 1 &  \\
\rowcolor{black!80!brown!15} 44 & 112 & 93 & 93 & 40 & 1 & 7 & 7 & 4 &  \\
\rowcolor{black!80!brown!15} 45 & 286 & 270 & 270 & 93 & 5 & 10 & 8 & 2  & \\
\rowcolor{black!80!brown!15} 46 & 15 & 13 & 13 & 7 & 1 & 5 & 5 & 2 &  \\
\rowcolor{black!80!brown!15} 47 & 4 & 4 & 4 & 2 & 4 & 2 & 2 & 1 &  \\
\rowcolor{black!80!brown!15} 48 & 15305 & 7890 & 7330 & 1394 & 1 & 12 & 13 & 52  & \\
\bottomrule
\end{tabular}

\caption{Schurian association schemes and their properties.}
\label{table:AllValues}
\end{table}
}
Up to order $30$, the values in the total, stratifiable, symmetric, and primitive columns agree with the values given in Table 2 of \cite{CAP} where they correspond to the number of $2$-closed permutation groups which are transitive, stratifiable, generously transitive, and primitive, respectively.
 The highlighted rows indicate that the enumeration of the corresponding Schurian association schemes is new, and is not available at \cite{hanaki_miyamoto}. The number of non-Schurian association schemes is given, as found at \cite{hanaki_miyamoto}, in column \emph{NS} for completeness. They indicate that the total number of association schemes may be significantly larger than the number which are Schurian, particularly as the order grows. As a result, further enumeration of association schemes is likely to be increasingly difficult. For example, there are $32730$ strongly regular graphs (which are equivalent to $2$-class association schemes) with $36$ vertices \cite{SRGs}, already far exceeding the total number of Schurian association schemes of this order.

Having computed the $2$-closed groups we were able to further compute the $2$-closure (up to conjugacy in $S_n$) of every transitive group of degree $n$, for $2 \leqslant n \leqslant 48$. This was achieved by utilising the following observation: if $G$ and $G'$ are groups such that $G'$ is $2$-closed, then $G'$ is conjugate in $S_n$ to the $2$-closure of $G$ if and only if $\mathcal{K}(G) \cong \mathcal{K}(G')$. Identifying the $2$-closure of a group $G$ is then a matter of testing isomorphism of $\mathcal{K}(G)$ against the elements of $L$. These isomorphism checks were also performed using \emph{AssociationSchemes}. Recall that an isomorphism between two association schemes $(\Omega, \mathcal{R})$ and $(\Delta, \mathcal{S})$ with $d$ classes is given by bijections $f:\Omega \to \Delta$ and $g: \{0, \ldots d\} \to \{0, \ldots d\}$ such that $(x, y) \in R_i$ if and only if $(f(x), f(y)) \in S_{g(i)}$. In \emph{AssociationSchemes} an isomorphism is given by $[\sigma_f, \sigma_g]$ where $\sigma_f$ and $\sigma_g$ are permutations corresponding to $f$ and $g$. Indeed, if $[\sigma_f, \sigma_g]$ is the isomorphism from $K(G')$ to $K(G)$ then the automorphism group of $K(G')$, and hence the $2$-closure of $G'$, is given by $G^{\sigma_f}$.

The results of this paper exploit the speed at which \emph{AssociationSchemes} is able to compute the automorphism group of an association scheme. Finding the $2$-closure of a related group is typically much slower and more difficult. We achieve this by representing an association scheme as an edge-coloured digraph where $(a, b)$ is a directed edge with colour i if and only if $(a, b) \in R_i$. Such a digraph is complete and satisfies the following property: the number of coloured triangles on a given directed edge depends only on the choice of colours. The automorphism group of the association scheme is simply the automorphism group of the corresponding directed graph. Following \cite[\S 14]{nug286}, an edge-coloured digraph may be represented as a vertex-coloured graph with ${\rm log}(d+1)$ layers. The automorphism group of the original digraph is then the found by taking the action of the automorphism group of this layered digraph on its first layer. This method is used by \emph{AssociationSchemes} to compute the automorphism group of an association scheme. Since it calls \emph{bliss} \cite{bliss} or \emph{nauty} \cite{nauty} to find the automorphism group, it is able to do this very quickly in practice. This method avoids computing intersections of groups or stabilisers of sets and so avoids many significant computational bottlenecks.

\emph{AssociationSchemes} is also able to compute canonical forms and isomorphisms of association schemes by considering isomorphisms of the corresponding graphs. A key difficulty in this process is accounting for the edge colourings (the order of the relations), since these may be relabelled to give isomorphic association schemes without the corresponding edge-coloured graphs being isomorphic. Prior to Version 3.0.0 of \emph{AssociationSchemes}, we overcame this by first computing all the algebraic automorphisms (reordering of relations preserving the intersection numbers) and applying the ismorphism/canonisation process to each of the edge-coloured digraphs that resulted. From Version 3.0.0, we instead allow colours to be exchanged within the digraph, combining the ideas in \cite[\S 14]{nug286} for edge-colouring and interchangable vertex-colouring. This requires a digraph with $d$ layers to be constructed, with an additional vertex added for each colour present. Calling \emph{bliss} \cite{bliss} on this digraph establishes an ordering on the relations. Accounting for this ordering on relations, a second call to \emph{bliss} then determines an ordering on the vertices. The result is a fast, practical method for computing canonical forms of association schemes and isomorphism between association schemes.

\section{Database}

The complete list of Schurian association schemes of order $2 \leq n \leq 48$ is made available as a supplementary database at \cite{SmallSchurianSchemesDatabase}. The database consists of files with the name \emph{SchurianSchemesN} where $N$ is the order of the association schemes. Each line of each file has the form
\[
[M, S, x, T, L],
\]
where:
    \begin{itemize}
\item $M$ is the relation matrix of the association scheme.
\item $S$ is a list of generators (as permutations) for the automorphism group of the association scheme; note that if 
$G = \langle S \rangle$ then $\mathcal{K}(G)$ is the association scheme defined by
$M$.
\item The integer $x$ is the transitive identification of the automorphism group $G$ in the \emph{Magma} and \emph{GAP} libraries\footnote{The transitive groups use the same ordering in both Magma and GAP, so the transitive identification numbers match.}.
\item $T$ is the character table of the association scheme. Note that some cyclotomic numbers are described as sums of roots of unity with rational coefficients, where the primitive $n$-th root of unity $e^{2 \pi/n}$ is written in the \emph{GAP} notation $E(n)$.
\item The $i$-th entry of
$L$ is the multiplicity of the character corresponding to the $i$-th row of 
$T$.
\end{itemize}

Consider, for example, the $3$-class association scheme with order $8$ defined by the relation matrix
\[ M =
\begin{bmatrix}
  0 & 1 & 2 & 2 & 3 & 2 & 1 & 1  \\
  2 & 0 & 2 & 1 & 1 & 3 & 1 & 2  \\
  1 & 1 & 0 & 1 & 2 & 2 & 3 & 2  \\
  1 & 2 & 2 & 0 & 2 & 1 & 1 & 3  \\
  3 & 2 & 1 & 1 & 0 & 1 & 2 & 2  \\
  1 & 3 & 1 & 2 & 2 & 0 & 2 & 1  \\
  2 & 2 & 3 & 2 & 1 & 1 & 0 & 1  \\
  2 & 1 & 1 & 3 & 1 & 2 & 2 & 0  
\end{bmatrix}.
\]
It is the Schurian scheme $\mathcal{K}(G)$, where
\[
G = \langle (1,3,5,7)(2,4,6,8), (1,3,8)(4,5,7) \rangle,
\]
with transitive identification $12$.
Its character table is
\[
T = 
\begin{bmatrix}
  1 & 3 & 3 & 1  \\
  1 & \sqrt{3}i & -\sqrt{3}i & -1  \\
  1 &-\sqrt{3}i & \sqrt{3}i & -1  \\
  1 & -1 & -1 & 1 
\end{bmatrix},	
\]
with multiplicities $[1, 2, 2, 3]$.
Observe, $\sqrt{3}i = e^{2\pi/3} - e^{4\pi/3}$ which is E(3) - E(3)\string^2 in \emph{GAP} notation. 
The corresponding entry in the database can be found on line $12$ of the file ``SchurianSchemes8'' where it appears as:\newline

\noindent{\tt
[ [ [ 0, 1, 2, 2, 3, 2, 1, 1 ], [ 2, 0, 2, 1, 1, 3, 1, 2 ], [ 1, 1, 0, 1, 2, 2, 3, 2 ], [ 1, 2, 2, 0, 2, 1, 1, 3 ], 
 [ 3, 2, 1, 1, 0, 1, 2, 2 ], [ 1, 3, 1, 2, 2, 0, 2, 1 ], [ 2, 2, 3, 2, 1, 1, 0, 1 ], [ 2, 1, 1, 3, 1, 2, 2, 0 ] ], \newline
 [ (1,3,5,7)(2,4,6,8), (1,3,8)(4,5,7) ], 12, [ [ 1, 3, 3, 1 ], [ 1, E(3)-E(3)\string^2, -E(3)+
  E(3)\string^2, -1 ], [ 1, -E(3)+E(3)\string^2, E(3)-E(3)\string^2, -1 ], [ 1, -1, -1, 1 ] ], [ 1, 2, 2, 3 ] ].
}\newline

The properties given in Table \ref{table:AllValues} are not given in the database because they are easily computable from the data provided. The character tables are included in the database, however, since they are incredibly useful in applications but can be very slow to compute. Note that it is significantly easier to verify a character table than to compute it.

The $2$-closures of the transitive permutation groups of degree $2 \leq n \leq 48$ are made available as a supplementary database at \cite{TwoClosuresDatabase}. The database consists of files with the name \emph{TwoClosuresN} where $N$ is the degree of the permutation groups.  Each line of each file has the form
\[
[x,y],
\]
where:
\begin{itemize}
\item The integer $x$ is the transitive identification of the group.
\item The integer $y$ is the transitive identification of the group isomorphic to the $2$-closure of the group with transitive identification $x$.
\end{itemize}
For  example, line $5491$ of the file ``TwoClosures48'' appears as:\newline\newline
\noindent{\tt
[5491,271829]
}\newline\newline
This says that the $2$-closure of the group accessed by {\tt TransitiveGroup(48, 5491)} is conjugate in $S_{48}$ to the group accessed by {\tt TransitiveGroup(48, 271829)}. Note that a group is $2$-closed if $x=y$. This data is likely to be particularly useful in the case where $n=48$, since there is no means of determining the transitive identification of a group in the \emph{GAP} or \emph{Magma} libraries for this degree.

\bibliographystyle{abbrv}
\bibliography{references.bib}

\end{document}